\documentclass [a4paper]{article}

\usepackage{amsfonts}
\usepackage{graphicx}
\usepackage{amsmath}
\usepackage{amssymb}
\usepackage{epsfig}

\newtheorem{thm}{Theorem}
\newtheorem{propn}[thm]{Proposition}
\newtheorem{lemma}{Lemma}

\newtheorem{cor}{Corollary}

\newcommand{\abs}[1]{\left\vert#1\right\vert}
\newcommand{\set}[1]{\left\{#1\right\}}

\newcommand{\Naturals}{\mathbb N}

\newcommand{\Area}{\text{Area}}
\newcommand{\Rad}{\text{Rad}}
\newcommand{\Star}{\text{star}}
\newcommand{\Diam}{\text{Diam}}

\begin{document}

\title{Filling Length in Finitely Presentable Groups}
\author{S. Gersten \and T. Riley}
\date{8 June 2000}

\maketitle

\begin{abstract}
Filling length measures the length of the contracting closed loops
in a null-homotopy.  The filling length function of Gromov for a
finitely presented group measures the filling length as a function
of length of edge-loops in the Cayley 2-complex. We give a bound
on the filling length function in terms of the log of an
isoperimetric function multiplied by a (simultaneously realisable)
isodiametric function.
\end{abstract}

\section{Isoperimetric and isodiametric functions}
Given a finitely presented group $\Gamma= \left<\mathcal{A} \left|
\mathcal{R} \right.\right>$ various \emph{filling invariants}
arise from considering reduced words $w$ in the free group
$F(\mathcal{A})$ such that $w=_{\Gamma}1$.  Such
\emph{null-homotopic} words are characterised by the existence of
an equality in  $F(\mathcal{A})$
\begin{equation}\label{eqn}
w=\prod_{i=1}^{N}{u_i^{-1} r_i u_i}
\end{equation}
for some $N \in \Naturals$, relators $r_i \in R$, and words $u_i
\in F(\mathcal{A})$.  A \emph{van Kampen diagram} provides a
geometric means of displaying such an equality - see \cite[page
155]{BrH}, \cite[pages 235ff]{Lyndon}.  This gives a notion of a
\emph{homotopy disc} for $w$.  Then, in analogy with
null-homotopic loops in a Riemannian manifold, we can associate
various filling invariants to the possible van Kampen diagrams for
null-homotopic words $w$.

Many such filling invariants are discussed by Gromov in Chapter 5
of \cite{Gromov}.  We will be concerned with three: the \emph{Dehn
function} (also known as the \emph{optimal isoperimetric
function}), the \emph{optimal isodiametric function} and the
\emph{filling length function}.

The first two of these invariants are better known than the third
- see for example \cite{Gersten2} and \cite{Gersten}. Let $|w|$
denote the length of a reduced word in $F(\mathcal{A})$. If
$w=_{\Gamma}1$ then define $\Area(w)$ to be the minimum number $N$
such that there is a van Kampen diagram for $w$ with $N$ faces
(i.e. 2-cells). The diameter of a van Kampen diagram $\mathcal{D}$
is the supremum over all vertices $v$ of $\mathcal{D}$ of the
shortest path in the 1-skeleton of $\mathcal{D}$ that connects $v$
to the base point of $\mathcal{D}$.  Define $\Diam(w)$ to be the
least diameter of van Kampen diagrams for $w$.  Then the Dehn
function $f_0:\Naturals \to \Naturals$ and optimal isodiametric
function $g_0:\Naturals \to \Naturals$ for $\Gamma=
\left<\mathcal{A} \left| \mathcal{R} \right.\right>$ are defined
by
\begin{eqnarray*}
f_0(n)&:=&\max\set{\Area(w) : w \in F(\mathcal{A}), |w| \leq n
\text{ and } w=_{\Gamma}1},\\ g_0(n)&:=&\max\set{\Diam(w) : w \in
F(\mathcal{A}), |w| \leq n \text{ and } w=_{\Gamma}1}.
\end{eqnarray*}
We say that $f$ and $g$ are respectively isoperimetric and
isodiametric functions for $\Gamma$ if $f_0(n) \leq f(n)$ and
$g_0(n) \leq g(n)$ for all $n$.

As defined the Dehn and the optimal isodiametric function are
dependent on the choice of presentation of $\Gamma$.  However
different presentations produce $\simeq$-equivalent\footnote{Given
two functions $f_{1},f_{2}:(0,\infty) \to (0,\infty)$ we say $
f_{1}\preceq f_{2}$ when there exists $C>0$ such that for all $l
\in (0,\infty)$, $f_{1}(l)\leq Cf_{2}(Cl+C)+Cl+C$. This yields the
equivalence relation: $f_{1} \simeq f_{2}$ if and only if $f_{1}
\preceq f_{2}$ and $f_{2} \preceq f_{1} $. } functions.  From an
algebraic point of view the Dehn function $f_0(n)$ is the least
$N$ such that for any null-homotopic word $w$ with $|w| \leq n$
there is an equality in $F(\mathcal{A})$ of the form of equation
(\ref{eqn}).  Similarly (but this time only up to
$\simeq$-equivalence) the optimal isodiametric function is the
optimal bound on the length of the conjugating words $u_i \in
F(\mathcal{A})$.

\section{The filling length function}\label{The filling length function}
In the context of a Riemannian manifold $X$ consider contracting a
null-homotopic loop $\gamma: [0,1] \to X$ based at $x_0 \in X$ to
the constant loop at $x_0$. By definition there is some continuous
$H:[0,1] \times [0,1] \to X$ denoted by $H_t(s)=H(t,s)$ with
$H_0=\gamma$, $H_1(s) = x_0$ for all $s$, and $H_t(0) = H_t(1) =
x_0$ for all $t$. \emph{Filling length} is a control on the length
of the loops $H_t$.  So the filling length of $H$ is the supremum
of the lengths of the loops $H_t$ for $t \in [0,1]$.  And the
filling length of $\gamma$ is the infimum of the filling lengths
of all possible null-homotopies $H$.  So (using Gromov's notation)
define $\text{Fill}_0 \text{ Leng} \,\ell$ to be the supremum of
the filling lengths of all null-homotopic loops $\gamma: [0,1] \to
X$ of length at most $\ell$ and based at $x_0 \in X$.

Now translate to the situation in a finitely presented group
$\Gamma= \left<\mathcal{A} \left| \mathcal{R} \right.\right>$. The
concept of homotopy discs is provided by van Kampen diagrams. In
the combinatorial context of a van Kampen diagram $\mathcal{D}$
 we use \emph{elementary homotopies}.  The boundary word $w$ is reduced to
the constant word $1$ at the base point $x_0$ by successively
applying two types of moves:
\begin{enumerate}\label{homotopy moves}
    \item \emph{(1-cell collapse)} remove pairs $(e^1,e^0)$ for which $e^0 \in \partial e^1$ is a
    0-cell which is not the base point $x_0$,
    and $e^1$ is a 1-cell only attached to the rest of the
    diagram at one 0-cell which is not $e^0$;
    \item \emph{(2-cell collapse)} remove pairs $(e^2,e^1)$
    where $e^2$ is a 2-cell of $\mathcal{D}$ with $e^1$ an edge of $\partial e^2 \cap \partial
    \mathcal{D}$ (note this does not change the 0-skeleton of $\mathcal{D}$).
\end{enumerate}
Algebraically 1-cell collapse corresponds to free reduction in
$F(A)$ (but not \emph{cyclic} reduction since the base point $x_0$
is preserved). And 2-cell collapse is the substitution of $a_2 a_3
\ldots a_n$ for $a_1^{-1}$, where $a_1, a_2, \ldots, a_n \in
\mathcal{A}^{\pm 1}$ and $a_1 a_2 \ldots a_n$ is a cyclic
permutation of an element of $\mathcal{R}^{\pm 1}$. Let the
filling length of $\mathcal{D}$, denoted $FL(\mathcal{D})$, be the
best possible bound on the length of the boundary word as we
successively apply these two types of move to reduce $w$ to $1$.
Define the filling length $FL(w)$ of a null-homotopic word $w$ by
$$FL(w):=\min \set {FL(\mathcal{D}) : \mathcal{D} \text{ is a van
Kampen diagram for } w }.$$ Then define the filling length
function $h_0:\Naturals \to \Naturals$ by $$h_0(n):=\max\set{FL(w)
: w \in F(\mathcal{A}), |w| \leq n \text{ and } w=_{\Gamma}1}.$$

As for $f_0$ and $g_0$, observe that $h_0$ is independent of the
presentation up to $\simeq$-equivalence.

\medskip
Some relationships known between $f_0,g_0$ and $h_0$ are as follows.

\medskip

\noindent{\textbf{Examples}
\begin{enumerate}
\item{\label{easy inequality}For a finitely presented group $\Gamma= \left<\mathcal{A}
\left| \mathcal{R} \right.\right>$ with $K:=\max \set{|r| : r \in
\mathcal{R}}$, we see that for all $n$ $$g_0(n) \leq h_0(n) \leq
2Kf_0(n)+n.$$ The first inequality arises since the concentric
loops of length at most $h_0(n)$ can be followed to reach the base
point.  The second inequality holds because given a
null-homotopic word $w$ of
length at most $n$, we find $2Kf_0(n)+n$ is at least twice the
total length of the 1-skeleton of a van Kampen diagram for $w$.}
\item{Filling length has also been used by Gersten and Gromov (see
pages 100ff of \cite{Gromov}) to give an isoperimetric function:
$$f_0 \preceq \exp h_0.$$  A
null-homotopic word $w$ of length $n$ can be reduced
to the identity by elementary homotopies through distinct words of
length at most $h_0(n)$. There are at most $\exp (C\,h_0(n))$ such
words for some constant $C>0$. This bounds the number of
\emph{2-cell collapse} moves and so this number is at least
$f_0(n)$.}
\item{Cohen and Gersten have given a double exponential bound
for $f_0$ in terms of $g_0$ \label{double exponential bound}(see
\cite{Cohen} and \cite{Gersten2}): $$ f_0 (n)\preceq \exp \exp
(g_0(n)+n).$$}
\item \emph{Asynchronously combable groups}
\label{Asynchronously combable groups} have linear bounds on their
filling length functions. This is a result of Gersten
\cite[Theorem 3.1 on page 130]{Gersten3} where the notation
$LCNH_1$ is in this case what we call linearly bounded filling
length. In essence the homotopy can be performed by contracting in
the direction of the combing, so the contracting loop always
remains normal to the combing lines. (See also \cite{Gersten} for
definitions.)
\end{enumerate}

\section{Logarithmic shelling of finite rooted trees}\label{shelling}

We now digress to a lemma about rooted trees. Let $\mathcal{T}$ be
a finite rooted tree in which each node has valence three except
for the root (valence two) and the leaves (valence one).

Let $\mathcal{F}$ be a finite forest of such trees. The
\emph{visible} nodes of $\mathcal{F}$ are the roots. An
\emph{elementary shelling} is the removal of the root of one of
its trees (together with the two edges that meet that root when
the tree has more than one node). A \emph{(complete) shelling} is
a sequence of elementary shellings ending with the empty forest.
The \emph{visibility number} of a shelling of $\mathcal{F}$ is the
maximum number of visible vertices occurring in the shelling. The
visibility number $VN(\mathcal{F})$ is the minimum visibility
number of all shellings.

Let $N(\mathcal{T})$ denote the number of nodes of $\mathcal{T}$.

\begin{lemma} Let the integer $d$ be determined by $2^d-1<N(\mathcal
T)\le 2^{d+1}-1$. Then $VN(\mathcal T)\leq d+1$.
\end{lemma}

\noindent \emph{Proof.} To obtain this bound on $VN(\mathcal T)$
we shall perform each elementary shelling by always choosing a
tree with the least number of nodes to shell first.

We argue by induction on $N(\mathcal T)$, where the induction
begins when $N(\mathcal T)=1$; in this case $d=0$ and $VN(\mathcal
T)=1$, as required.

For the induction step, assume that $N(\mathcal T)>1$ with
$2^d-1<N(\mathcal T)\leq 2^{d+1}-1$. Removing the root of
$\mathcal T$ produces two trees $\mathcal T_1, \mathcal T_2$. We
let $N(\mathcal T_1)\leq N(\mathcal T_2)$. Let
$2^{d_i}-1<N(\mathcal T_i)\leq 2^{d_i+1}-1$ for $i=1,2$. By the
induction hypothesis we have $VN(\mathcal T_i)\leq d_i+1$ for
$i=1,2$. Since we shell $T_1$ first, we get $VN(\mathcal T)\leq
\max(VN(\mathcal T_1)+1,VN(\mathcal T_2)) \leq \max(d_1+2,d_2+1)$
by the induction hypothesis. There are now two cases, depending on
whether $d_1<d_2$ or $d_1=d_2$.

\medskip
\noindent \emph{Case 1}. $d_1<d_2$. In this case $\max
(d_1+2,d_2+1)=d_2+1\le d+1$, so we get $VN(\mathcal T)\leq d+1$ as
required.

\medskip
\noindent \emph{Case 2}. $d_1=d_2$.  Here
$\max(d_1+2,d_2+1)=d_1+2$. We have $ 2(2^{d_1}-1)+1<N(\mathcal
T_1)+N(\mathcal T_2)+1\leq 2(2^{d_1+1}-1)+1$, whence
$2^{d_1+1}-1<N(\mathcal T)\leq 2^{d_1+2}-1$. It follows that
$d=d_1+1$, and $VN(\mathcal T)\leq d+1$ as required.

\medskip
This completes the induction, and the proof of Lemma 1 is
complete.

\begin{cor}\label{corol}  $VN(\mathcal T)< \log_2( N(\mathcal T) +
1) + 1$.
\end{cor}

\noindent \emph{Proof.}  Write $2^d-1<N(\mathcal T)\leq
2^{d+1}-1$, so $VN(\mathcal T)\leq d+1< \log_2(N(\mathcal
T)+1)+1$, as required.

\medskip
\noindent \emph{Remark.} Note that since $VN(\mathcal T)$ is an
integer, the upper bound for $VN(\mathcal T)$ in the corollary can
be replaced by the least integer bounding $\log_2( N(\mathcal T) +
1)$ from above. Stated in this form, the result is sharp, as we
see by taking $T$ to be the complete rooted tree $T(d)$ of depth
$d$.  In this case $T(d)$ has $2^{d+1}-1$ nodes, and the
visibility number is $d+1$.

\section{A bound on the filling length function}\label{main thm
section}

We now proceed towards our main theorem.
\medskip

\noindent \textbf{Definition}  Let $\mathcal P$ be a finite
presentation for the group~$\Gamma$. An AD-pair for $\mathcal P$
is a pair of functions $(f,g)$ from $\Naturals$ to $\Naturals$
such that for every circuit $w$ of length at most $n$ in the
Cayley graph there exists a van Kampen diagram $\mathcal{D}_w$
with area at most $f(n)$ and diameter at most $g(n)$. Note that
$f$ is an isoperimetric function and $g$ is an isodiametric
function.

\medskip
\noindent{\textbf{Examples} Up to common multiplicative constants
the following are examples of AD-pairs for groups $\Gamma$.}

\begin{enumerate}

\item  $(E^x, x)$ for some $E>1$ and $(x L(x),x)$ are both AD-pairs
when $\Gamma$ is an asynchronously combable group.  Here the
length function $L(n)$ is the maximum length of combing paths
 for group elements at distance at most $n$ from
the identity. That $(E^x,x)$ is an AD-pair follows from the linear
bound on the filling length function and that $f_0 \preceq \exp
h_0$; see section 2 Example 4.
In particular $(x^2, x)$ is an AD-pair when $\Gamma$ is
synchronously automatic since $\Gamma$ then admits a combing in
which the combing lines are quasi-geodesics (see \cite[pages
84-86]{Epstein}).

\item  $(f(x),f(x)+x)$ where $\Gamma$ is an arbitrary
 finitely presented group with an isoperimetric
function $f$. (See Gersten \cite{Gersten}, Lemma 2.2.)

\item $(E^{E^{g(x)+x}},g(x))$ for some $E>1$, when $\Gamma$ is an
arbitrary finitely presented group with an isodiametric function $g$.
(See Gersten \cite{Gersten2}.)

\item $(x^r,x^{r-1})$ where $\Gamma$ satisfies a polynomial
isoperimetric function of degree $r \geq 2$.  We postpone proof of
this example to section 5.

\item \label{Bridson's groups}$(x^{2m+1},x^m)$ where $\Gamma_m$ is
\begin{eqnarray*}
\lefteqn{ \left<a_1,\ldots,a_m,s,t,\tau \left|\,\text{for } i <
m,\, s^{-1} a_i s=a_{i+1},\right.\right.} \\ & &
\hspace{1in}\left.
[t,a_i]=[\tau,a_i]=[s,a_m]=[t,a_m]=[\tau,a_mt]=1 \right>.
\end{eqnarray*} This family of examples is due to Bridson - see
\cite{Bridson}.  He shows that in fact $x^{2m+1}$ is the
\emph{optimal} isoperimetric function and $x^m$ is the
\emph{optimal} isodiametric function.

\end{enumerate}

\medskip

\begin{thm} \label{Main Theorem}
Let $(f,g)$ be an AD-pair for the finite presentation $\mathcal
P$. Then $h_0(n) \preceq g(n)\log (f(n)+1)$ for all $n$.
\end{thm}

We actually prove the stronger statement:
\begin{propn} \label{Proposition}
Suppose $\Gamma$ is a finitely and triangularly presented group,
and $w \in \Gamma$ is null-homotopic with $n:=|w|$. Given a van
Kampen diagram $\mathcal{D}$ for $w$ with $D:=\Diam(\mathcal{D})$
and $A:=\Area(\mathcal{D})$ we find $$FL(\mathcal{D}) \leq
(2D+1)(\log_2(A+1)+1)+4D+1+n.$$
\end{propn}

Any finite presentation $\left<\mathcal{A} \left| \mathcal{R}
\right.\right>$ for a group $\Gamma$ yields a finite
\emph{triangular} \label{triangular presentation} presentation for
$\Gamma$. Such presentations are characterised by the length of
relators being at most three. If $r \in \mathcal{R}$ is
expressible in $F(\mathcal{A})$ as $w_1w_2$ where $|w_1|,|w_2|
\geq 2$ then add a new generator $a$ to $\mathcal{A}$, and in
$\mathcal{R}$ replace $r$ by $a^{-1}w_1$ and $aw_2$.  A triangular
presentation is achieved after a finite number of such
transformations.

As $f_0,g_0$ and $h_0$ are invariant up to $\simeq$-equivalence on
change of finite presentation, Proposition \ref{Proposition} is
sufficient to prove Theorem \ref{Main Theorem}.

\medskip
\noindent \emph{Proof of Proposition \ref{Proposition}.} We start
by taking a maximal geodesic tree $\mathcal{T}$ in the 1-skeleton
of the van Kampen diagram $\mathcal{D}$, and rooted at the base
point $x_0$ of $\mathcal{D}$. So from any vertex of $\mathcal{D}$
there is a path in $\mathcal{T}$ to $x_0$ with length at most $D$.

By cutting along paths in $\mathcal{T}$ we can decompose
$\mathcal{D}$ into sub-diagrams $\mathcal{D}_i$ where only one
edge from $\partial\mathcal{D} - \mathcal{T}$ occurs in each
$\mathcal{D}_i$.  We will perform the elementary homotopy which
realises the filling length bound by pushing across each of these
$\mathcal{D}_i$ in turn, as shown in Figure \ref{fig:Elementary
homotopy}.  Then
\begin{equation}\label{eqn2}
\text{FL}(\mathcal{D}) \leq \displaystyle{\max_i}
\set{\text{FL}(\mathcal{D}_i)} + n.
\end{equation}

\begin{figure}[ht]
\centerline{\epsfig{file=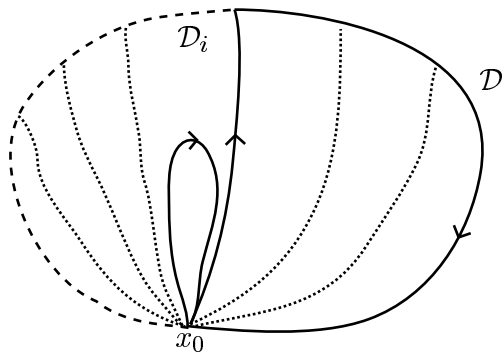}}
\caption{Elementary homotopy of
$\mathcal{D}$.}\label{fig:Elementary homotopy}
\end{figure}

It remains to explain how to perform the elementary homotopy
across each $\mathcal{D}_i$.  We will use six types of
\emph{2-cell collapse} moves (see
section 2 above).
These are depicted in Figure 2, with solid lines representing
edges in $\mathcal{T}$.

\begin{figure}[ht]
\centerline{\epsfig{file=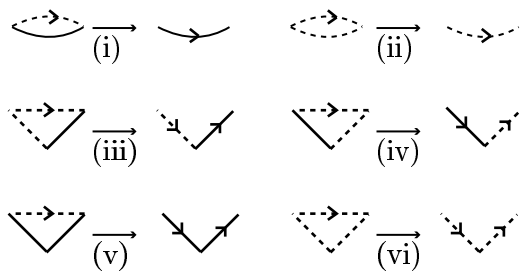}}
\caption{Homotopy
moves.}\label{fig:Homotopy moves}
\end{figure}

We now describe the means of performing the homotopy in a way that
controls filling length. Repeatedly apply the following four
steps:
\begin{enumerate}
    \item 1-cell collapse (see section 2 above),
    \item moves (i) and (ii): bi-gon collapse,
    \item moves (iii) and (iv),
    \item moves (v) and (vi) in accordance with \emph{logarithmic
    shelling}.
\end{enumerate}
The first step in the list that is available is performed, and
then we return to the start of the list.  Repeating this, the
boundary loop of $\mathcal{D}_i$ will eventually be reduced to the
constant loop at the base point $x_0$.

The means by which we use logarithmic shelling to choose which
2-cell to collapse when performing step 4 requires some
explanation.  The result of cutting $\mathcal{D}_i$ along
$\mathcal{T}$ and removing 2-cells of the type encountered in
(i),(ii),(iii), and (iv) of Figure \ref{fig:Homotopy moves}, is
illustrated in Figure \ref{fig:Shelling tree}. Taking the dual
gives a rooted tree of the form discussed in Section
\ref{shelling}.  The 2-cell to be pushed across is chosen in
accordance with the process of logarithmic shelling of rooted
trees discussed there.  The number of nodes in the tree is at most
$A$ and so by Corollary \ref{corol} the visibility number is at
most $\log_2 (A+1)+1$.

\begin{figure}[ht]
\centerline{\epsfig{file=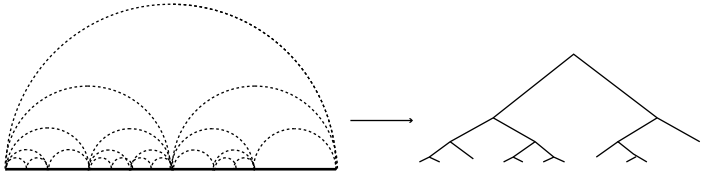}}
\caption{The shelling
tree of $\mathcal{D}_i$.}\label{fig:Shelling tree}
\end{figure}

It remains to explain how performing the elementary homotopy as
described above leads to the required bound on $h_0$. Consider the
situation when the next step to be applied is number 4. The
visibility number associated to the dual tree described above is
at most $\log_2 (A+1)+1$.  The homotopy loop includes at most
$\log_2 (A+1)+1$ edges of the type occurring in move (vi). These
are separated by paths in $\mathcal{T}$ of length at most $2D$.
The loop is closed by another path in $\mathcal{T}$ again of
length at most $2D$.  So this loop has length at most:
\[\log_2 (A+1)+1+2D\log_2 (A+1)+2D=(2D+1)(\log_2(A+1)+1).\]
Now applying move (v) or (vi) increases the length of the loop by
1, creating two new \emph{channels} where moves (i), (ii), (iii)
and (iv) may be performed. Consider then applying steps 1, 2 and
3. Step 1 can only decrease the length, and step 2 leaves it
unchanged. Step 3 can be applied at most $2D$ times in each of the
two channels.  Thus the increase in length before step 4 is next
applied is at most $1+4D$.  This gives bound
\[\text{FL}(\mathcal{D}_i)\leq(2D+1)(\log_2(A+1)+1)+4D+1.\]
Combining this with the inequality (\ref{eqn2}) we have our
result.

\section{Examples}\label{Examples section}
In this section we provide a proof for the assertion made in
section \ref{The filling length function} about an AD-pair for a
group with a polynomial isoperimetric function. We then give
applications of Theorem \ref{Main Theorem} to particular classes
of groups, and we conclude with a discussion of an open question.

\medskip
\begin{propn} \label{AD-pair proposition} Let $\Gamma$ be a group admitting a polynomial
isoperimetric function of degree $r \geq 2$.  Then up to a common
multiplicative constant $(x^r,x^{r-1})$ is an AD-pair for
$\Gamma$.
\end{propn}

Papasoglu gives this result for $r=2$ in \cite[page
799]{Papasoglu}. It requires a small generalisation of his
argument to obtain the result for all $r \geq 2$, as follows. (See
also \cite[page 100]{Gromov}.)

We will make use of some of definitions. Let the radius of a van
Kampen diagram $\mathcal{D}$ to be
$$\Rad(\mathcal{D}):=\max\set{d(v,\partial \mathcal{D}):v \text{
is a vertex of }\mathcal{D}},$$ where $d(v,\partial \mathcal{D})$
is the combinatorial distance in the 1-skeleton.

For a subcomplex $K$ of $\mathcal{D}$ define $\Star(K)$ to be the
union of closed 2-cells meeting $K$.  Define $\Star_i(K)$ to be
the $i$-th iterate of the star operation for $i\ge 1$;
by convention $\Star_0(K)=K$. So if $\Gamma$ is
triangularly presented then the 0-cells in $\Star_i(\partial
\mathcal{D})$ are precisely those a distance at most $i$ from
$\partial D$.

\medskip
The substance of Proposition \ref{AD-pair proposition} is in the
following lemma.

\begin{lemma}\label{AD-pair
lemma} Suppose $\Gamma= \left<\mathcal{A} \left| \mathcal{R}
\right.\right>$ is triangularly presented (see
the paragraph following Proposition 2)
and that $\mathcal{R}$ includes
all null-homotopic words of length at most 3. Suppose further that
there is $M>0$ such that $f_0(n) \leq Mn^r$ for all $n$. Then for
all null-homotopic $w$ we have $\Rad(\mathcal{D})\leq
12M\abs{w}^{r-1}$, where $\mathcal{D}$ is a minimal area van
Kampen diagram for $w$.
\end{lemma}

As discussed in section 4,
a change
of finite presentation induces a $\simeq$-equivalence on
isoperimetric and isodiametric functions. Note that there are only
finitely many words of length at most 3 in a finitely presented
group.  Also observe that adding $n/2$ is sufficient to obtain a
diameter bound from a radius bound. Thus this lemma is sufficient
to prove Proposition \ref{AD-pair proposition}.

\medskip
\noindent \emph{Proof of Lemma \ref{AD-pair lemma}.} We proceed by
induction on $n$. For $n \leq 3$ the result follows from our
insistence that $\mathcal{R}$ includes all null-homotopic words of
length at most 3.

For the induction step suppose $w$ is null-homotopic and $\abs{w}
= n$. Let $\mathcal{D}$ be a minimal area van Kampen diagram for
$w$. Let $N_i:=\Star_i(\partial \mathcal{D})$. Let $c_i:=\partial
N_i - \partial \mathcal{D}$, which is the union of simple closed
curves any two of which meet at one point or not at all (note that
$c_0= \partial(\Star_0(\partial D)) -\partial D$, which is empty).

Now $\Area(N_{i+1})-\Area(N_i) \geq l(c_i)/3$ because every 1-cell
of $c_i$ lies in the boundary of some 2-cell in $N_{i+1}-N_i$.

For all $i$, $$\Area(\mathcal{D}) \geq \Area(N_{i+1}) \geq
l(c_i)/3 + l(c_{i-1})/3+ \cdots + l(c_0)/3. $$  Thus if $l(c_i) >
n/2$ for all $i\leq 6Mn^{r-1}$ we get a contradiction of the area
bound $Mn^r$ for $\mathcal{D}$.  So for some $i \leq 6Mn^{r-1}$ we
find $l(c_i) \leq n/2$.  We can appeal to the inductive hypothesis
to learn that the diagrams enclosed by the simple closed curves
constituting $c_i$ have radius at most $12M(n/2)^{r-1}\leq
6Mn^{r-1}$.

So a vertex $v$ of $\mathcal{D}$ either lies in $N_i$, in which
case $d(v,\partial \mathcal{D}) \leq 6Mn^{r-1}$, or is in a
diagram enclosed by one of the simple closed curves $c$ of $c_i$ .
In the latter case $d(v,\partial \mathcal{D}) \leq d(v, c)+
d(c,\partial \mathcal{D}) \leq 12Mn^{r-1}$ as required,
thus completing the
proof of the lemma.

\medskip
We now give some applications of our main theorem to
particular classes of groups.

\medskip
\noindent \textbf{Example} \emph{Polynomial isoperimetric
function.}  If the finitely presented group $\Gamma$ admits a
polynomial isoperimetric function of degree $r \geq 2$, it follows
from Theorem \ref{Main Theorem} and Proposition \ref{AD-pair
proposition} above that $h_0(n) \preceq n^{r-1} \log(n+1)$.  This
contrasts with the inequality $h_0 \preceq f_0$ in section 2
Example 1.

\medskip
\noindent \textbf{Example} \emph{Bridson's Groups $\Gamma_m$} (see
section 3 Example 5).
 These have
AD-pairs $(x^{2m+1},x^m)$ and so Theorem \ref{Main Theorem} gives
us bounds of $x^m\log x$ on their filling length functions, which
is a significant improvement on the bounds $x^{2m+1}$ obtained
from the inequality $h_0 \preceq f_0$.

\medskip
\noindent \textbf{Open question}  In connection with the double
exponential bound quoted in section 2 Example 3,
 it is an open problem, to our knowledge first raised by
John Stallings, whether there is always a simple exponential bound
$f_0 \preceq \exp g_0$.  It is natural then to ask whether
there is always an AD-pair of the form $(\exp g_0, g_0)$.
(This adds the requirement that the $\exp g_0$ bound on $f_0$ is always
realisable on the same van Kampen diagram as $g_0$.) Our main
theorem gives a necessary condition that this be true, namely that
$h_0 \preceq {g_0}^2$.

\medskip

We shall now make some observations relevant to the single
exponential question just stated.

\begin{propn} If $\mathcal{P}$ is a finite presentation,
then for all integers $N\ge 3$ there exists $C(N)>0$
such that for all van Kampen diagrams $\mathcal{D}$ in
$\mathcal{P}$ all of whose vertices have valence at most
$N$ one has
\begin{enumerate}
\item $Area(\mathcal{D})\le N^{Diam(\mathcal{D})+1}-1$, and
\item $FL(\mathcal{D})\le C(N)\cdot Diam(\mathcal{D})^2+n+1$.
\end{enumerate}
\end{propn}

\noindent \emph{Proof.} The number of vertices at a given distance
$i$ from the base point is at most $N(N-1)^{i-1}$, so it follows
that the number of geometric edges $E(\mathcal{D})$ satisfies
$E(\mathcal{D})\le N+N^2+\cdots +N^D < \frac{N^{D+1}-1}{N-1}$,
where $D=Diam(\mathcal{D})$.  Since each edge is incident with at
most 2 faces, we get $Area(\mathcal{D})\le
2\frac{N^{D+1}-1}{N-1}\le N^{D+1}-1$, giving the first conclusion
of the proposition.

>From Proposition 2 it follows that
$FL(\mathcal{D})\le (2D+1)(D+1)\log_2(N)+4D+1+n
=C(N)D^2+n+1$,
where $C(N)$ depends only on $N$, proving the second conclusion.

\medskip
\begin{cor}
For every finite presentation $\mathcal{P}$
there is a constant $C>0$
such that if $\mathcal{D}$ is an immersed topological
disc diagram in $ \mathcal{P}$, then
$FL(\mathcal{D})\le C\cdot\Diam(\mathcal{D})^2+n+1$.
\end{cor}

\noindent \emph{Proof.}
Since $\mathcal{D}$ is immersed, the valence of a vertex $v$ is
at most the number of
edges incident at a vertex of the Cayley graph,
namely, twice the number of generators of
$\mathcal{P}$.  The corollary follows from the
second conclusion of the proposition.

\medskip
\noindent \emph{Remark.}  Gromov observed in \cite{Gromov} 5C
that if $h_0\preceq g_0$, then
it follows that $f_0\preceq \exp g_0$; one sees this
as a consequence of section 2 Example 3 above.
We do not know an example from finitely presented groups
where $h_0\preceq g_0$ fails;
however it is shown in \cite{Frankel-Katz} that this
can fail in a simply connected Riemannian
context.
\footnote{Their example does not amount to a properly discontinuous
cocompact action by isometries on a simply connected Riemannian
manifold, so it does not correspond to an example arising from finitely
presented groups.}

\bibliographystyle{plain}
\bibliography{bibli}

\begin{thebibliography}{10}

\bibitem{Bridson}
M.~Bridson.
\newblock Asymptotic cones and polynomial isoperimetric inequalities.
\newblock {\em Topology}, 38(3):543--554, 1999.

\bibitem{BrH}
M.~Bridson and A.~Haefliger.
\newblock {\em Metric Spaces of Non-positive Curvature}.
\newblock Number 319 in Grundlehren der mathematischen Wissenschaften. Springer
  Verlag, 1999.

\bibitem{Cohen}
D.~E. Cohen.
\newblock Isoperimetric and isodiametric inequalities for group presentations.
\newblock {\em Int. J. of Alg. and Comp.}, 1(3):315--320, 1991.

\bibitem{Epstein}
D.~B.~A. Epstein, J.~W. Cannon, S.~F. Holt, S.~V.~F. Levy, M.~S. Paterson, and
  W.~P. Thurston.
\newblock {\em Word Processing in Groups}.
\newblock Jones and Bartlett, 1992.

\bibitem{Frankel-Katz}
S.~Frankel and M.~Katz.
\newblock The {M}orse landscape of a {R}iemannian disc.
\newblock {\em Ann. Inst. Fourier, Grenoble}, 43(2):503--507, 1993.

\bibitem{Gersten2}
S.~Gersten.
\newblock The double exponential theorem for isoperimetric and isodiametric
  functions.
\newblock {\em Int. J. of Alg. and Comp.}, 1(3):321--327, 1991.

\bibitem{Gersten}
S.~Gersten.
\newblock Isoperimetric and isodiametric functions.
\newblock In G.~Niblo and M.~Roller, editors, {\em Geometric group theory II},
  number 182 in LMS lecture notes. Camb. Univ. Press, 1993.

\bibitem{Gersten3}
S.~Gersten.
\newblock Asynchronously automatic groups.
\newblock In Charney, Davis, and Shapiro, editors, {\em Geometric group
  theory}, pages 121--133. de Gruyter, 1995.

\bibitem{Gromov}
M.~Gromov.
\newblock Asymptotic invariants of infinite groups.
\newblock In G.~Niblo and M.~Roller, editors, {\em Geometric group theory II},
  number 182 in LMS lecture notes. Camb. Univ. Press, 1993.

\bibitem{Lyndon}
R.~C. Lyndon and P.~E. Schupp.
\newblock {\em Combinatorial Group Theory}.
\newblock Springer Verlag, 1977.

\bibitem{Papasoglu}
P.~Papasoglu.
\newblock On the asymptotic invariants of groups satisfying a quadratic
  isoperimetric inequality.
\newblock {\em J. Diff. Geom.}, 44:789--806, 1996.

\end{thebibliography}

\begin{tabbing}
E-mail: gersten@math.utah.eduuuuuu\= E-mail: rileyt@maths.ox.ac.uk
\kill  \textsc{Steve Gersten} \> \textsc{Tim Riley}\\ \small
Mathematics Department \> \small Mathematical Institute\\ \small
University of Utah\> \small 24-29 St Giles \\ \small Salt Lake
City \> \small Oxford
\\ \small UT 84112 \> \small OX1 3LB\\ \small USA \> \small UK
\\ \small \hspace{3mm} E-mail: gersten@math.utah.edu \> \small \hspace{3mm} E-mail: rileyt@maths.ox.ac.uk
\end{tabbing}

\end{document}